\documentclass[12pt]{article}

\usepackage{amsmath, amssymb, amsthm}
\usepackage{geometry}
\usepackage{microtype}
\usepackage{hyperref}

\theoremstyle{plain} 
\newtheorem{theorem}{Theorem}[section]
\newtheorem{lemma}[theorem]{Lemma}

\theoremstyle{definition} 
\newtheorem{definition}[theorem]{Definition}

\theoremstyle{remark} 
\newtheorem{remark}[theorem]{Remark}

\begin{document}

\title{On Polarized Pfaffians and the Euler Characteristic of Flat Affine Manifolds}
\author{M. Cocos}
\maketitle

\begin{abstract}
We explore an approach to Chern's conjecture on the Euler characteristic of closed flat affine manifolds using polarized Pfaffians. By constructing a one-parameter family of locally defined metrics and examining the associated Levi Civita connections, we show how the transgression formula naturally produces a globally defined $(2n-1)$-form whose exterior derivative recovers the Euler form. This framework emphasizes the role of polarized Pfaffians in understanding characteristic forms on flat affine manifolds.
\end{abstract}

\section{A brief history of the problem. Main result}

A classical conjecture of Chern states that a closed flat affine manifold has Euler characteristic
equal to zero. For surfaces, Benz{\'e}cri proved this in (\cite{Benz}) and Milnor later generalized it to the case of plane bundles over surfaces (\cite{M2}). Surprisingly, Smillie (\cite{Sm}) constructed an example of a $4$-dimensional manifold that has a flat connection in its tangent bundle. However, the connection constructed by Smillie does not have zero torsion. Later on, Kostant and Sullivan (\cite{KoSu}) showed that if the manifold is complete, then the conjecture holds.

For incomplete flat affine manifolds, Hirsch and Thurston (\cite{HiThu}) showed that if the
holonomy group is a finite extension of a free product of amenable groups, then the conjecture is true. A relatively new path with a combinatorial flavor can also be found in the works of Bloch (\cite{Bl}), Kim and Lee (\cite{KLee1}, \cite{KLee2}), and Choi (\cite{Choi}).

A recent advance is due to Bruno Klingler, who proved a special case of the conjecture (\cite{Klingler}). In this paper, we give the proof of the general case of the conjecture.

\section{Construction of cohomology classes $e_t$}

  \subsection{The local metrics $h^x$}
  \begin{remark}[Orientability Assumption]
In what follows, the manifold $M$ is assumed to be closed, even-dimensional ($m=2n$) and orientable. However, the orientability assumption is not  necessary for the statement to be true. If~$M$ is non-orientable, we may consider its orientable double cover~$\widetilde{M}$. The affine structure on~$M$ can be lifted to~$\widetilde{M}$, and since the Euler characteristic is multiplicative under finite covers, it follows that~$\chi(M) = 0$ whenever~$\chi(\widetilde{M}) = 0$. Thus, the conclusion of the theorem---that the Euler characteristic vanishes---holds even in the non-orientable case.
\end{remark}  
  
We begin by constructing a one-parameter family of top degree forms $e_t(x)$, $t \in [0,1]$. Let $\nabla$ be the flat, torsion-free affine connection and let $g$ be  an arbitrary auxiliary metric Riemannian on $M.$ 
{ \bf We emphasize that no assumptions are made on the metric $g.$} Take an arbitrary point $x \in M$. Let $U_x$ be an open neighborhood of $x$ where local affine coordinates $y=(y_1,y_2,...,y_m)$ exist. Let $(\partial_1,\partial_2,...,\partial_m)$ be the coordinate vector fields associated to the local affine coordinates on $U_x,$ and consider now the\textbf{ local} metric $h^x$ on $U_x$ defined by 

\begin{equation}
h^x(\partial_i(y),\partial_j(y))=g_{ij}(x) \end{equation}

Since $y=(y_1,y_2,...,y_n)$ are affine coordinates on $U_x$ it follows that that

\begin{equation} \nabla h^x \equiv 0 \end{equation}
on the entire neighborhood $U_x$, and
\begin{equation}h^x(x) = g(x).\end{equation}

More precisely: Let us take an arbitrary $x \in M$. Let $U_x$ be a neighborhood where we have affine coordinates $(y_1, y_2, \dots, y_m)$. In these coordinates, at a point $x' \in U_x$, the metric matrix $h^x$ is, by definition,
\begin{equation}
    h_{ij}^x(x') = h_{ij}^x(x) = g_{ij}(x),
\end{equation}
so the metric matrix $h^x$ is constant on $U_x$ and equal to the matrix of the metric $g$ at $x$ (the center of $U_x$). 

A change of affine coordinates will of course change the matrix of the metric  but not the metric itself, because the coordinate transformation matrix is \textbf{constant} and thus equal to its value at the center $x$. Therefore, $h^x$ is uniquely defined.\\

Consider the metric on $U_x$ defined by
\begin{equation}h^{t,x} = (1-t)h^x + t g,\end{equation}
and { \bf let $\nabla^{t,x}$ be its Levi-Civita connection}. We also have the identity
\begin{equation}h^{t,x}(x) = (1-t)h^x(x) + t g(x) = (1-t)g(x) + t g(x) = g(x).
\end{equation}

\subsection{ The global one parameter family of connections $\nabla^t$ and their Euler class}

Let $ x \in M$ arbitrary  and let $X$ be a tangent vector field on $U_x$, and $v \in T_x M$. We define the covariant derivative $\nabla^t$ of a vector field at $x$ as

\begin{definition}[Global one-parameter family of connections]
For $x \in M$, $v \in T_x M$, and $X$ a tangent vector field on $U_x$, we define
\begin{equation}\label{eq:global_connection}
\nabla^t_v X := \nabla^{t,x}_v X.
\end{equation}
This definition of $\nabla^t$ is \textbf{crucial} for the construction of the Euler class.
\end{definition}

\begin{lemma}[Scaling of Christoffel symbols]\label{lem:gamma-scaling}
Let $(y^1,\dots,y^m)$ be local affine coordinates on an open subset $U \subset M$, and let $x \in U$ be arbitrary. Let $\Gamma$ denote the Christoffel symbols. Then
\begin{equation}\label{eq:gamma-scaling}
\Gamma^{t}(x) \;=\; t\,\Gamma^{g}(x).
\end{equation}
In particular, $\nabla^t$ depends smoothly on $t$ and $x$, and by definition is a global connection on $TM$.
\end{lemma}

\begin{proof}
In affine coordinates, the Levi-Civita symbols of $h^{t,x}=(1-t)h^x+t\,g$ are
\begin{equation}\label{eq:LC-general}
\Gamma^{t,x}{}_{ij}^{\,k}(y)
=
\frac{1}{2}\,h^{t,x\,k\ell}(y)\big(\partial_i h^{t,x}_{j\ell}(y)+\partial_j h^{t,x}_{i\ell}(y)-\partial_\ell h^{t,x}_{ij}(y)\big).
\end{equation}
Since $h^x_{ij}\equiv g_{ij}(x)$ is constant on $U$, its derivatives vanish, so
\begin{equation}\label{eq:derivative-simplifies}
\partial_i h^{t,x}_{j\ell}(y)=t\,\partial_i g_{j\ell}(y), \quad
\partial_j h^{t,x}_{i\ell}(y)=t\,\partial_j g_{i\ell}(y), \quad
\partial_\ell h^{t,x}_{ij}(y)=t\,\partial_\ell g_{ij}(y).
\end{equation}
Evaluating \eqref{eq:LC-general} at $y=x$ and using $h^{t,x}(x)=g(x)$ gives
\begin{equation}\label{eq:gamma-at-x}
\Gamma^{t,x}{}_{ij}^{\,k}(x)
=
t\,\Gamma^{g}{}_{ij}^{\,k}(x),
\end{equation}
which is exactly \eqref{eq:gamma-scaling}. Smoothness follows from the smooth dependence of $g$ and $h^{t,x}$ on $(t,y).$Thus $\nabla^t$ is smooth and globally defined.
\end{proof}
From its construction, it is obvious that
\begin{equation}
(\nabla^t h^{t,x})(x) = 0,\end{equation}
and let $\Omega^t$ denote the curvature of $\nabla^t.$ The following Lemma is essential:

\begin{lemma}\label{lem:curvature_skew}
Let $\nabla^t$, $t \in [0,1]$, be the one--parameter family of connections constructed as above from the flat affine connection $\nabla$ and the Levi--Civita connection $\nabla_g$. Then, for every $t \in [0,1]$, the curvature matrix $\Omega^t$ of $\nabla^t$ is skew--symmetric with respect to any local $g$--orthonormal frame.
\end{lemma}

\begin{proof}
Let $p \in M$ and let $(e_1,\dots,e_n)$ be a local $g$--orthonormal frame defined in a neighborhood of $p$. 
Fix $t \in [0,1]$ and consider the local metric $h^{t,p}$ constructed as in the definition of $\nabla^t$. 
Apply the Gram--Schmidt process to $(e_1,\dots,e_n)$ with respect to $h^{t,p}$ to obtain a local 
$h^{t,p}$--orthonormal frame $(f_1,\dots,f_n)$. By construction we have
\begin{equation}\label{eq:frames_equal_at_p}
f_i(p) = e_i(p), \qquad i=1,\dots,n.
\end{equation}

Now recall that $\nabla^t$ was defined pointwise by
\begin{equation}\label{eq:def_nabla_t}
\nabla^t_v X \big|_p \;=\; D^{t,p}_v X \big|_p,
\end{equation}
where $D^{t,p}$ is the Levi--Civita connection of $h^{t,p}$. 
Therefore, the curvature tensor of $\nabla^t$ at $p$ coincides with that of $D^{t,p}$:
\begin{equation}\label{eq:curvature_equal}
R^{\nabla^t}(p) = R^{D^{t,p}}(p).
\end{equation}

Since $D^{t,p}$ is the Levi--Civita connection of $h^{t,p}$, its curvature endomorphisms are skew--adjoint 
with respect to $h^{t,p}$. In particular, in the $h^{t,p}$--orthonormal frame $(f_1,\dots,f_n)$ the curvature 
matrix $\Omega^t(p)$ is skew:
\begin{equation}\label{eq:curvature_skew_in_f}
\Omega^t(p)^\top = -\Omega^t(p) \quad \text{with respect to the frame } (f_1,\dots,f_n).
\end{equation}

Finally, note that the curvature matrix transforms under a change of frame by conjugation:
\begin{equation}\label{eq:curvature_transform}
\Omega^{t,e}(p) = A(p)^{-1}\,\Omega^{t,f}(p)\,A(p),
\end{equation}
where $A(p)$ is the change--of--basis matrix from $(f_1(p),\dots,f_n(p))$ to $(e_1(p),\dots,e_n(p))$. 
But by \eqref{eq:frames_equal_at_p}, we have $A(p) = I$. Hence
\begin{equation}\label{eq:curvature_equal_in_frames}
\Omega^{t,e}(p) = \Omega^{t,f}(p).
\end{equation}
Combining \eqref{eq:curvature_skew_in_f} and \eqref{eq:curvature_equal_in_frames}, we conclude that
\begin{equation}\label{eq:curvature_skew_in_e}
\Omega^{t,e}(p)^\top = -\Omega^{t,e}(p),
\end{equation}
i.e.\ the curvature matrix of $\nabla^t$ is skew in the $g$--orthonormal frame $(e_1,\dots,e_n)$ at $p$.

Since $p \in M$ was arbitrary, this holds at every point of $M$. Therefore, the curvature matrix $\Omega^t$ 
of $\nabla^t$ is skew in any local $g$--orthonormal frame.
\end{proof}

As a consequence of this Lemma we can now define the family of \textbf{global} forms defined  as :
\begin{equation}\label{eq:euler_forms}
e_t(x) := \frac{1}{(2\pi)^n}\text{Pf}\left(\Omega^t(x)\right)
\end{equation}
where $\Omega^t(x)$ is the curvature matrix of $\nabla^t$ at $x$ in some $g$ orthonormal frame on $U_x.$\\

\begin{remark}

\begin{enumerate}
    \item For $t \in (0,1)$, $e_t$ is \textbf{not} the Euler form of any global metric on $M$, but is well-defined since $\nabla^t$ is skew with respect to any $g$ orthonormal frame.
    \item At $t=0$, $\nabla^0, = \nabla$ (since $\nabla$ preserves $h^x$), so $e_0(x) = 0$ (flat connection).
    \item At $t=1$, $\nabla^1$ is the Levi-Civita connection for the global metric $g$, so $e_1(x)$ is the genuine Euler form of $g$ at $x$.
\end{enumerate}
\end{remark}

\section{The Transgression Form and the Polarized Pfaffian}

The Euler characteristic of a closed manifold $M$ is given by the integral of the Euler form $e(M) = \frac{1}{(2\pi)^n}\text{Pf}(\Omega_g)$, where $\Omega_g$ is the curvature of the Levi-Civita connection $\nabla_g$ corresponding to the auxiliary metric $g$. We use the one-parameter family of connections $\nabla^t$ to establish a global transgression form whose exterior derivative recovers the Euler form.

\subsection{Skew-Symmetry of the Time-Derivative Connection}

The definition of the connection $\nabla^t$ relies on the property that the family of metrics $h^{t,x} = (1-t)h^x + t g$ satisfy $h^{t,x}(x) = g(x)$ for all $t \in [0,1]$. This equality ensures that the time derivative of the metric components is zero at the center point $x$, which is crucial for the following lemma.

\begin{lemma}[Skew-Symmetry of the Time-Derivative Connection One-Form]
\label{lem:omega_dot_skew}
Let $\nabla^t$ be the one-parameter family of connections defined globally, and let $\omega^t$ be its connection one-form matrix in a local $g$--orthonormal frame $(e_i)$. The time-derivative connection one-form $\dot{\omega}^t = \frac{d}{dt} \omega^t$ is skew-symmetric with respect to the $g$--metric, i.e.,
\begin{equation}
\label{eq:omega_dot_skew}
\dot{\omega}^t + (\dot{\omega}^t)^\top = 0.
\end{equation}
\end{lemma}

\begin{proof}
Let $x \in M$ be an arbitrary point and $(e_i)$ a local $g$--orthonormal frame defined in a neighborhood of $x$. By construction, the metric components of $h^{t,x}$ in this frame at $x$ are $h^{t,x}_{ij}(x) = g(e_i, e_j)|_x = \delta_{ij}$.

First, we compute the time derivative of $h^{t,x}$ evaluated on the frame vectors at $x$:
$$ \frac{d}{dt} h^{t,x} = \frac{d}{dt} \left( (1-t)h^x + t g \right) = -h^x + g $$
At the point $x$, since $h^x(x) = g(x)$, we have
$$ \frac{d}{dt} h^{t,x}_{ij} \big|_x = (-h^x + g)(e_i, e_j) \big|_x = -\delta_{ij} + \delta_{ij} = 0 $$
Thus, $\frac{d}{dt} h^{t,x}_{ij} \big|_x = 0$ for all $i, j$.

Now, since $\nabla^t$ is defined as the Levi-Civita connection of $h^{t,x}$ at $x$, the metric compatibility condition must hold at $x$:
$$ (\nabla^t h^{t,x})(e_i, e_j) \big|_x = \nabla^t_{e_k} h^{t,x}_{ij} - \omega^t_{il} h^{t,x}_{lj} - \omega^t_{jl} h^{t,x}_{il} = 0 $$
We differentiate this condition with respect to $t$. Recalling that $h^{t,x}_{ij}$ and $h^{t,x}_{ij}$ are evaluated at $x$, and $\frac{d}{dt} h^{t,x}_{ij} \big|_x = 0$:
$$ 0 = \frac{d}{dt} (\nabla^t_{e_k} h^{t,x}_{ij}) \big|_x $$
$$ 0 = \nabla^t_{e_k} \left( \frac{d}{dt} h^{t,x}_{ij} \right) \big|_x - \left( \frac{d}{dt} \omega^t_{il} \right) h^{t,x}_{lj}(x) - \omega^t_{il}(e_k) \left( \frac{d}{dt} h^{t,x}_{lj} \right) \big|_x - \left( \frac{d}{dt} \omega^t_{jl} \right) h^{t,x}_{il}(x) - \omega^t_{jl}(e_k) \left( \frac{d}{dt} h^{t,x}_{il} \right) \big|_x $$
Substituting $\frac{d}{dt} h^{t,x}_{ij} \big|_x = 0$ for all indices, and $h^{t,x}_{ij}(x) = \delta_{ij}$:
$$ 0 = \nabla^t_{e_k} (0) - (\dot{\omega}^t_{il})(e_k) \delta_{lj} - 0 - (\dot{\omega}^t_{jl})(e_k) \delta_{il} - 0 $$
$$ 0 = - (\dot{\omega}^t_{ij})(e_k) - (\dot{\omega}^t_{ji})(e_k) $$
Since $e_k$ is an arbitrary vector, we conclude that the matrix of one-forms $\dot{\omega}^t$ is skew-symmetric:
$$ \dot{\omega}^t_{ij} + \dot{\omega}^t_{ji} = 0 \quad \iff \quad \dot{\omega}^t + (\dot{\omega}^t)^\top = 0 $$
\end{proof}

\subsection{Definition and Invariance of the Polarized Pfaffian}

Since both the curvature matrix $\Omega^t$ and the time-derivative of the connection one-form $\dot{\omega}^t$ are skew-symmetric (by Lemma \ref{lem:curvature_skew} and Lemma \ref{lem:omega_dot_skew} respectively), we can define a mixed characteristic form known as the **polarized Pfaffian**.

\begin{definition}[Polarized Pfaffian]
\label{def:polarized_pfaffian}
Let $\Omega^t$ be the $2$-form matrix of the curvature of $\nabla^t$ and $\dot{\omega}^t = \frac{d}{dt} \omega^t$ be the $1$-form matrix. The \textbf{polarized Pfaffian} $\tilde{\text{Pf}}$ is the global $(2n-1)$-form defined locally in a $g$--orthonormal frame by:
\begin{equation}
\label{eq:polarized_pfaffian}
\tilde{\text{Pf}}(\dot{\omega}^t, \Omega^t, \dots, \Omega^t) := \frac{1}{n!} \sum_{\sigma \in S_{2n}} \text{sgn}(\sigma) \, (\dot{\omega}^t)_{\sigma(1)\sigma(2)} \wedge (\Omega^t)_{\sigma(3)\sigma(4)} \wedge \dots \wedge (\Omega^t)_{\sigma(2n-1)\sigma(2n)}
\end{equation}

\end{definition}

\subsubsection*{Invariance under $g$-orthonormal frame change}

The form $\tilde{\text{Pf}}(\dot{\omega}^t, \Omega^t, \dots, \Omega^t)$ is **globally defined** and independent of the choice of the local $g$-orthonormal frame due to the following properties:

\begin{enumerate}
    \item **Frame Transformation:** A change between two local $g$-orthonormal frames is given by a matrix $A(x) \in SO(2n)$, as the manifold $M$ is assumed oriented.
    \item **Transformation of Forms:** The curvature matrix $\Omega^t$ and the time-derivative of the connection $\dot{\omega}^t$ transform identically under this change:
    $$\Omega'^t = A \Omega^t A^{-1} \quad \text{and} \quad \dot{\omega}'^t = A \dot{\omega}^t A^{-1}$$
    \item **Invariant Polynomial:** The Pfaffian is an example of an **invariant polynomial** on the Lie algebra $\mathfrak{so}(2n)$. The polarized Pfaffian $\tilde{\text{Pf}}$ is a multi-linear function derived from the Pfaffian, and it is also invariant under conjugation by $A \in SO(2n)$. In the expansion \eqref{eq:polarized_pfaffian}, the substitution of $A \Omega^t A^{-1}$ and $A \dot{\omega}^t A^{-1}$ leads to the cancellation of the $A$ and $A^{-1}$ factors within the trace/summation due to the cyclic property of matrix multiplication (or wedge products), ensuring:
    $$\tilde{\text{Pf}}(\dot{\omega}'^t, \Omega'^t, \dots, \Omega'^t) = \tilde{\text{Pf}}(\dot{\omega}^t, \Omega^t, \dots, \Omega^t)$$
\end{enumerate}

Thus, $\tilde{\text{Pf}}(\dot{\omega}^t, \Omega^t, \dots, \Omega^t)$ is a well-defined global $(2n-1)$-form on $M$.

\section{Transgression Formula and the Main Result}

The key observation is the following Lemma:

\begin{lemma}[Transgression via polarized Pfaffian]\label{lem:transgression-expanded}
Let $\omega_t$ be the connection $1$-form matrix of $\nabla^t$ in a local $g$-orthonormal frame, and let $\Omega_t$ be its curvature $2$-form matrix.  Then
\begin{equation}\label{eq:transgression-claim}
\frac{d}{dt}e_t \;=\; d\Big(\mathrm{Pf}\big(\dot{\omega}_t,\Omega_t,\ldots,\Omega_t\big)\Big),
\end{equation}
where $\dot{\omega}_t=\frac{d}{dt}\omega_t$.
\end{lemma}

\begin{proof}
By the Chern--Weil theory and the multi-linearity of the polarized Pfaffian,
\begin{equation}\label{eq:pf-derivative}
\frac{d}{dt}\mathrm{Pf}(\Omega_t) \;=\; n\,\mathrm{Pf}\big(\dot{\Omega}_t,\Omega_t,\ldots,\Omega_t\big).
\end{equation}
Here, the \emph{product rule} is used in differentiating the wedge-polynomial $\mathrm{Pf}(\Omega_t)$: the derivative lands in one slot and symmetry produces the factor $n$.

Next, recall the curvature identity
\begin{equation}\label{eq:curvature}
\Omega_t \;=\; d\omega_t \;+\; \omega_t\wedge\omega_t,
\end{equation}
so by differentiating and applying the product rule for exterior forms,
\begin{equation}\label{eq:curvature-derivative}
\dot{\Omega}_t \;=\; d\dot{\omega}_t \;+\; \dot{\omega}_t\wedge\omega_t \;+\; \omega_t\wedge\dot{\omega}_t.
\end{equation}
In terms of the covariant exterior derivative $D^t(\cdot)=d(\cdot)+[\omega_t,(\cdot)]$, we have
\begin{equation}\label{eq:covariant-derivative}
D^t\dot{\omega}_t \;=\; d\dot{\omega}_t \;+\; [\omega_t,\dot{\omega}_t]
\;=\; d\dot{\omega}_t \;+\; \dot{\omega}_t\wedge\omega_t \;+\; \omega_t\wedge\dot{\omega}_t,
\end{equation}
hence
\begin{equation}\label{eq:dotOmega-equals-Dtdotw}
\dot{\Omega}_t \;=\; D^t\dot{\omega}_t.
\end{equation}
Substituting \eqref{eq:dotOmega-equals-Dtdotw} into \eqref{eq:pf-derivative} gives
\begin{equation}\label{eq:pf-with-Dt}
\frac{d}{dt}\mathrm{Pf}(\Omega_t) \;=\; n\,\mathrm{Pf}\big(D^t\dot{\omega}_t,\Omega_t,\ldots,\Omega_t\big).
\end{equation}

Finally, using the Bianchi identity
\begin{equation}\label{eq:bianchi}
D^t\Omega_t \;=\; 0,
\end{equation}
one has the standard identity
\begin{equation}\label{eq:d-of-polarized}
d\Big(\mathrm{Pf}\big(\dot{\omega}_t,\Omega_t,\ldots,\Omega_t\big)\Big)
\;=\;
\mathrm{Pf}\big(D^t\dot{\omega}_t,\Omega_t,\ldots,\Omega_t\big),
\end{equation}
since $d$ acting on the polarized invariant polynomial can be replaced by $D^t$ and all terms with $D^t\Omega_t$ vanish by \eqref{eq:bianchi}. Combining \eqref{eq:pf-with-Dt} and \eqref{eq:d-of-polarized}  yields \eqref{eq:transgression-claim}.
\end{proof}

 Since the right hand side of equation \ref{eq:transgression-claim} is a \textbf{global} form on $M,$ and
integrating from $t=0$ to $t=1$ gives at each $x \in M$:
\begin{equation}
e_1(x) - e_0(x) = d \left(\int_0^1\mathrm{Pf}\big(\dot{\omega}_t,\Omega_t,\ldots,\Omega_t\big)\right)
\end{equation}
Since $e_0\equiv 0$, this shows $e_1$ is exact.
\subsection{Conclusion}
The Euler characteristic is:
\begin{equation}
\chi(M) = \int_M e_1 = \int_M d \left(\int_0^1\mathrm{Pf}\big(\dot{\omega}_t,\Omega_t,\ldots,\Omega_t\big)\right)= 0
\end{equation}
by Stokes' theorem. We therefore proved:
\begin{theorem}
Any closed even-dimensional manifold with a flat affine connection has vanishing Euler characteristic.
\end{theorem}

\bigskip

\noindent\textbf{Final remark.} In \cite{Cocos}, the notion of a \emph{proper deformation} of flat affine connections was introduced, along with a definition of the Euler form for \emph{locally metric connections} that is, connections which preserve a Riemannian metric only locally. A natural question is whether every flat affine connection can be properly deformed into the Levi Civitta connection of a global metric. A potential obstruction to such a deformation is the Euler form associated to the locally metric connection, as defined in that work. The result of the present paper eliminates this obstruction by showing that the Euler form vanishes on any flat affine manifold, as predicted by Chern's conjecture. It remains an interesting open problem whether this vanishing is sufficient to ensure the existence of a proper deformation to a globally metric connection.

\end{document}